\pdfoutput=1
\documentclass[11pt]{amsart}

\usepackage[utf8]{inputenc}
\usepackage[english]{babel}
\usepackage[a4paper,margin=2.3cm]{geometry}
\usepackage{amsmath,amssymb,amsthm,mathtools}
\usepackage{microtype}
\usepackage{needspace}
\usepackage[colorlinks=true,linkcolor=blue,citecolor=blue,urlcolor=black]{hyperref}
\hypersetup{pdftitle={Dirichlet deviations for Hill operators with two exponential terms of equal modulus},pdfauthor={Berkay Anahtarc\i}}

\newtheorem{theorem}{Theorem}[section]
\newtheorem{proposition}[theorem]{Proposition}
\newtheorem{lemma}[theorem]{Lemma}
\newtheorem{corollary}[theorem]{Corollary}
\theoremstyle{remark}
\newtheorem*{remark}{Remark}

\numberwithin{equation}{section}

\newcommand{\C}{\mathbb{C}}
\newcommand{\R}{\mathbb{R}}
\newcommand{\N}{\mathbb{N}}
\newcommand{\Z}{\mathbb{Z}}
\newcommand{\ls}{\lambda_n^{s}}
\newcommand{\lc}{\lambda_n^{c}}
\newcommand{\Th}{\Theta_n}
\newcommand{\Img}{\operatorname{Im}}

\title[Dirichlet deviations for Hill operators]{Dirichlet deviations for Hill operators with two exponential terms of equal modulus}
\author{Berkay Anahtarc{\i}}
\address{Department of Mathematical Engineering, Faculty of Engineering, \"{O}zye\u{g}in University, \.{I}stanbul, T\"{u}rkiye}
\email{berkay.anahtarci@ozyegin.edu.tr}

\begin{document}

\begin{abstract}
Consider the Hill operator $L=-d^2/dx^2+ae^{-2ix}+be^{2ix}$ on $[0,\pi]$ with $a,b\in\C\setminus\{0\}$. For large $n$, let $\mu_n$ be the Dirichlet eigenvalue and $\lambda_n^-$, $\lambda_n^+$ the two periodic ($n$~even) or antiperiodic ($n$~odd) eigenvalues near $n^2$, counted with algebraic multiplicity, and set $\delta_n=\mu_n-\frac12(\lambda_n^-+\lambda_n^+)$. For $|a|\neq|b|$, Djakov and Mityagin proved that $\delta_n\sim-2F_n(a^n+b^n)$, where $F_n=1/(4^n[(n-1)!]^2)$. With equal moduli, the two powers can nearly or exactly cancel, and their results give this formula along sequences of indices on which $(a/b)^n$ stays bounded away from $-1$. In this case we obtain the first correction to their formula and prove that, as $n\to\infty$ along all indices,
\[
\delta_n=-2F_n(a^n+b^n)+\frac{F_n}{n}(a-b)(a^n-b^n)+O\Bigl(\frac{F_n|a|^n}{n^2}\Bigr),
\]
with the error uniform in $a$ and $b$ for fixed $ab$. Consequently, their formula holds along all indices with $a^n+b^n\neq0$ if and only if $a/b$ is a root of unity. At the indices with $a^n+b^n=0$, the second term is the leading term of $\delta_n$, and in the non-self-adjoint case $b\neq\bar a$ the first two terms give the asymptotics of $\delta_n$ along all indices. We also locate $\mu_n$ relative to both edges up to $O(n^{-2})$ times the gap length, and show that the expansion is sharp for $a\neq b$.
\end{abstract}

\maketitle

\section{Introduction}\label{sec:intro}

Consider the Hill operator
\begin{equation}\label{eq:L}
L=-\frac{d^2}{dx^2}+v(x),\qquad v(x)=ae^{-2ix}+be^{2ix},\qquad a,b\in\C\setminus\{0\},
\end{equation}
on $[0,\pi]$. For large $n$, the Dirichlet problem for $L$ has exactly one eigenvalue $\mu_n$ near $n^2$, and the periodic ($n$~even) or antiperiodic ($n$~odd) problem has exactly two such eigenvalues, $\lambda_n^-$ and $\lambda_n^+$, counted with algebraic multiplicity \cite{DM20}. We call $\lambda_n^+-\lambda_n^-$ the $n$th spectral gap, $|\lambda_n^+-\lambda_n^-|$ its length, and $\lambda_n^\pm$ its edges. The eigenvalues $\lambda_n^-$, $\lambda_n^+$ and $\mu_n$ are the vertices of the $n$th spectral triangle. We locate $\mu_n$ by the edge deviations $\mu_n-\lambda_n^\pm$ and the midpoint deviation $\delta_n=\mu_n-\frac12(\lambda_n^-+\lambda_n^+)$. The latter does not depend on the labelling of $\lambda_n^\pm$.

For real periodic potentials, the periodic and antiperiodic problems are self-adjoint, and $\mu_n$ lies between $\lambda_n^-$ and $\lambda_n^+$, so the triangle reduces to a segment. In this case the decay of the gap lengths characterizes the smoothness and analyticity of the potential. For complex periodic potentials, such as \eqref{eq:L} with $b\neq\bar a$, these problems are not self-adjoint, and $\mu_n$ may lie off the segment joining the edges. The gaps alone then need not determine regularity \cite{Gas80}, and the position of the Dirichlet eigenvalues must be taken into account (see \cite{DM03,Tka92} and the references therein).

The deviations also govern the convergence of expansions in the eigenfunctions and associated functions of these two problems. In the non-self-adjoint case, the systems of these root functions need not contain a Riesz basis of $L^2([0,\pi])$, so such expansions may fail to converge unconditionally. By the criteria of Djakov and Mityagin \cite{DM12} and Gesztesy and Tkachenko \cite{GT12}, each of these systems contains such a basis exactly when the ratios $|\mu_n-\lambda_n^\pm|/|\lambda_n^+-\lambda_n^-|$ are bounded for large $n$ of the corresponding parity with $\lambda_n^+\neq\lambda_n^-$. Equivalently, $|\delta_n|/|\lambda_n^+-\lambda_n^-|$ is bounded, since $\mu_n-\lambda_n^\pm=\delta_n\mp\frac12(\lambda_n^+-\lambda_n^-)$. For the operator \eqref{eq:L}, this condition holds if and only if $|a|=|b|$ \cite{DM11}.

For $|a|\neq|b|$, Djakov and Mityagin \cite{DM20} showed that, as $n\to\infty$, the midpoint deviation satisfies
\begin{equation}\label{eq:DM}
\delta_n\sim-2F_n\,(a^n+b^n),
\end{equation}
where $F_n=1/(4^n[(n-1)!]^2)$. They write $\delta_n$ for $\mu_n-\lambda_n^+$. For all nonzero $a$ and $b$, Veliev's isospectrality result \cite{Vel13} and the Mathieu gap asymptotics \cite{AD12} give the gap only up to sign. Thus $(\lambda_n^+-\lambda_n^-)^2\sim64F_n^2(ab)^n$, and the gap length is asymptotic to $8F_n|ab|^{n/2}$. If $|a|\neq|b|$, then one of $a^n$ and $b^n$ dominates in \eqref{eq:DM}, so $|\delta_n|/|\lambda_n^+-\lambda_n^-|\to\infty$. Hence $\mu_n$ is far from both edges relative to the gap, and the Riesz-basis condition fails.

For $|a|=|b|$ the picture is different. The gap length is then asymptotic to $8F_n|a|^n$, and the modulus of the right-hand side of \eqref{eq:DM} is $\frac14|1+(a/b)^n|\,(1+o(1))$ times the gap length. Thus this term is comparable to the gap unless $(a/b)^n$ is near~$-1$, where it is small, and it vanishes when $(a/b)^n=-1$. Djakov and Mityagin obtain \eqref{eq:DM} along any sequence of indices on which $(a/b)^n$ stays bounded away from $-1$, and the leading terms of the edge deviations along any sequence on which it stays bounded away from $1$.

Their results locate $\mu_n$ in two limiting situations. When $(a/b)^n\to-1$, both edge formulas apply, and their average gives $\delta_n=o(|\lambda_n^+-\lambda_n^-|)$, so $\mu_n$ is close to the midpoint on the scale of the gap. When $(a/b)^n\to1$, \eqref{eq:DM} and the asymptotics of the squared gap give
\[
(\mu_n-\lambda_n^-)(\mu_n-\lambda_n^+)=\delta_n^2-\tfrac14(\lambda_n^+-\lambda_n^-)^2=o\bigl(|\lambda_n^+-\lambda_n^-|^2\bigr),
\]
so $\mu_n$ is close to an edge on the same scale.

If $a/b$ is a root of unity, then $(a/b)^n$ takes only finitely many values, so \eqref{eq:DM} holds along all indices except those of exact cancellation, where $a^n+b^n=0$ and the relation has no meaning. If $a/b$ is not a root of unity, then $a^n+b^n\neq0$ for all $n$, but $(a/b)^n$ comes arbitrarily close to $-1$ and to $1$.

These results leave three questions open. First, along which sequences of indices does \eqref{eq:DM} still hold? Second, what is the leading term of $\delta_n$ at exact cancellation, and what replaces \eqref{eq:DM} along sequences on which it fails? Third, how close does $\mu_n$ come to the nearer edge when $(a/b)^n$ is near $1$?

Our main result, Theorem~\ref{thm:main}, is that, for $|a|=|b|$,
\[
\delta_n=-2F_n(a^n+b^n)+\frac{F_n}{n}(a-b)(a^n-b^n)+O\bigl(F_n|a|^nn^{-2}\bigr)
\]
as $n\to\infty$ along all indices, with the error uniform in $a$ and $b$ for fixed $ab$. The second term is then the first correction to \eqref{eq:DM}. By Corollary~\ref{cor:sharp}, if $a\neq b$, the first term alone has an error of exact order $F_n|a|^n/n$ along infinitely many indices, so the correction cannot be omitted. The remainder of the two-term expansion is $o(F_n|a|^nn^{-2})$ if and only if $a=b$, and then it is asymptotic to $a^{n+2}F_n/n^3$. For fixed $a$ and $b$, along sequences of indices with $|(a/b)^n-1|=O(n^{-1})$, Theorem~\ref{thm:refined} gives the next term of the expansion, $-(a-b)^2a^nF_n/(2n^2)$, with remainder $O(F_n|a|^nn^{-3})$. The answers to the three questions are as follows.
\begin{enumerate}
\item[\textup{(i)}] Along a sequence of indices with $a^n+b^n\neq0$, \eqref{eq:DM} holds exactly when $n|1+(a/b)^n|\to\infty$, that is, when the distance of $(a/b)^n$ from $-1$ is much larger than $1/n$. It holds along all such indices exactly when $a/b$ is a root of unity (Theorem~\ref{thm:validity}).
\item[\textup{(ii)}] At exact cancellation, $\delta_n\sim2(a-b)a^nF_n/n$. If $b\neq\bar a$ or $a/b$ is a root of unity, the sum of the first two terms of the expansion approximates $\delta_n$ with relative error $O(n^{-1})$ along all indices (Proposition~\ref{prop:midpoint}).
\item[\textup{(iii)}] Let $a\neq b$. If $b\neq\bar a$ or $a/b$ is a root of unity, Proposition~\ref{prop:edges}(a) gives the edge deviations along all indices with relative error $O(n^{-1})$. At the indices with $a^n=b^n$, the leading expression of \cite{DM20} for the deviation from the nearer edge vanishes, and $\mu_n$ lies at a distance asymptotic to $|a-b|^2|\lambda_n^+-\lambda_n^-|/(16n^2)$ from that edge (Proposition~\ref{prop:edges}(b)).
\end{enumerate}
If $b=\bar a$ and $a/b$ is not a root of unity, the expansion keeps its absolute error bound, but we do not obtain a relative error along all indices. Finally, for all $a$ and $b$ with $|a|=|b|$, Proposition~\ref{prop:position} determines the ratios $|\mu_n-\lambda_n^\pm|/|\lambda_n^+-\lambda_n^-|$ of the Riesz-basis criteria up to $O(n^{-2})$.

The general asymptotics of the Dirichlet eigenvalues \cite{Mar11,KST12,Yil05} have remainder bounds that decay at most like a negative power of $n$, so these bounds cannot resolve the deviations at the factorial scale $F_n|a|^n$. However, $v$ is a real translate of a Mathieu potential if and only if $|a|=|b|$, and in that case a translation identity expresses the equation for $\mu_n$ through two functions of the Mathieu operator, one vanishing at each edge (Section~\ref{sec:shift}). Factoring out these zeros gives $\mu_n$ exactly as a weighted combination of the edges, so an error in the weight is multiplied by the gap, and the $O(n^{-2})$ solution estimates of Section~\ref{sec:asymptotics} suffice. Section~\ref{sec:gap} determines the sign of the gap when the edges are labelled by the parity of their Mathieu eigenfunctions. Section~\ref{sec:main} determines the weight and proves the main expansion, its refinement and its sharpness. Section~\ref{sec:consequences} derives the answers (i)--(iii) to the three questions above.

\section{Reduction to a translate of the Mathieu potential}\label{sec:shift}

Recall from \eqref{eq:L} that $v(x)=ae^{-2ix}+be^{2ix}$ with $a,b\in\C\setminus\{0\}$. From now on we assume $|a|=|b|$. Let $q\in\C$ satisfy $q^2=ab$. Then $|a/q|=1$, so there is $\tau\in[0,\pi)$ with $e^{2i\tau}=a/q$, and since $b=q^2/a$,
\begin{equation}\label{eq:shift}
a=qe^{2i\tau},\qquad b=qe^{-2i\tau}.
\end{equation}
With the Mathieu potential $u(x)=2q\cos2x$, \eqref{eq:shift} gives $v(x)=qe^{-2i(x-\tau)}+qe^{2i(x-\tau)}=u(x-\tau)$. For each $\lambda\in\C$, the translation $T\colon Y\mapsto Y(\cdot-\tau)$ therefore maps the solutions of the Mathieu equation $-y''+uy=\lambda y$ bijectively onto those of $-y''+vy=\lambda y$. It preserves $\pi$-periodicity and $\pi$-antiperiodicity, so $L$ and the Mathieu operator $-d^2/dx^2+u$ have the same periodic and antiperiodic eigenvalues. By contrast, it identifies the Dirichlet problem for $L$ with that for the Mathieu operator on $[-\tau,\pi-\tau]$, so for fixed $q$ only the Dirichlet spectrum can depend on $\tau$.

Since $q$ is determined up to sign, the pair $(q,\tau)$ is unique up to the change $(q,\tau)\mapsto(-q,\tau+\pi/2)$, with $\tau+\pi/2$ taken modulo~$\pi$. For later use we also note that $|q|=|a|=|b|$ and, for $n\in\N$,
\begin{equation}\label{eq:shift-cons}
a-b=2iq\sin2\tau,\qquad a^n=q^ne^{2in\tau},\qquad b^n=q^ne^{-2in\tau},\qquad (a/b)^n=e^{4in\tau}.
\end{equation}

Let $s_v(x,\lambda)$ solve $-y''+vy=\lambda y$ with $y(0)=0$, $y'(0)=1$, so that the Dirichlet eigenvalues of $L$ are the zeros of $s_v(\pi,\cdot)$. Let $M_v\colon y\mapsto y(\cdot+\pi)$ be the monodromy operator on the two-dimensional space of solutions of this equation, and let $\Delta_v(\lambda)=\operatorname{tr}M_v$ be the Hill discriminant. Since the Wronskian $W(y_1,y_2)=y_1y_2'-y_1'y_2$ of two solutions is constant, $\det M_v=1$, so $\pm1$ is an eigenvalue of $M_v$ exactly when $\Delta_v(\lambda)=\pm2$. Hence the periodic and antiperiodic eigenvalues are the zeros of $\Delta_v-2$ and $\Delta_v+2$, respectively.

To express $s_v(\pi,\cdot)$ and $\Delta_v$ through solutions of the Mathieu equation, let $C(x,\lambda)$ and $S(x,\lambda)$ be the solutions with initial values $C(0,\lambda)=S'(0,\lambda)=1$ and $C'(0,\lambda)=S(0,\lambda)=0$, so that $W(C,S)=1$. For each $x\in\R$, $s_v(x,\lambda)$, $C(x,\lambda)$, $S(x,\lambda)$ and their $x$-derivatives are entire functions of $\lambda$ by analytic dependence on parameters. Here and below, we suppress $\lambda$ from the notation when no confusion can arise. By uniqueness, every solution $y$ of the Mathieu equation equals $y(0)C+y'(0)S$. Applied to the solutions $C(\cdot+\pi)$ and $S(\cdot+\pi)$, this gives the action of the monodromy operator $M_u$ of $u$ on the basis $C,S$, namely, for all $x\in\R$,
\begin{equation}\label{eq:monodromy}
C(x+\pi)=C(\pi)C(x)+C'(\pi)S(x),\qquad S(x+\pi)=S(\pi)C(x)+S'(\pi)S(x).
\end{equation}
Since $u$ is even, $C(-\cdot)$ and $-S(-\cdot)$ are solutions with the same initial data as $C$ and $S$, so uniqueness gives
\begin{equation}\label{eq:parity}
C(-x)=C(x),\qquad S(-x)=-S(x).
\end{equation}

\needspace{5\baselineskip}
\begin{lemma}[monodromy identities]\label{lem:even}
For all $\lambda\in\C$, we have
\begin{enumerate}
\item[\textup{(a)}] $C(\pi,\lambda)=S'(\pi,\lambda)$ and $C(\pi,\lambda)^2-C'(\pi,\lambda)S(\pi,\lambda)=1$;
\item[\textup{(b)}] $\Delta_v(\lambda)=2C(\pi,\lambda)$ and $\Delta_v(\lambda)^2-4=4\,C'(\pi,\lambda)\,S(\pi,\lambda)$.
\end{enumerate}
\end{lemma}

\begin{proof}
(a) Since $u$ is even and $\pi$-periodic, $u(\pi-x)=u(x)$, so $S(\pi-\cdot)$ is also a solution. Hence the Wronskian $W(C,S(\pi-\cdot))=-C(x)S'(\pi-x)-C'(x)S(\pi-x)$ is constant. It equals $-S'(\pi)$ at $x=0$ and $-C(\pi)$ at $x=\pi$, so $C(\pi)=S'(\pi)$. With this, $W(C,S)=1$ at $x=\pi$ reads $C(\pi)^2-C'(\pi)S(\pi)=1$.

(b) For each solution $Y$ of the Mathieu equation, $M_v(TY)=Y(\cdot+\pi-\tau)=T(M_uY)$, so $M_v=TM_uT^{-1}$ and $\Delta_v=\operatorname{tr}M_u$. By \eqref{eq:monodromy}, the matrix of $M_u$ in the basis $C,S$ has the diagonal entries $C(\pi)$ and $S'(\pi)$, so (a)~gives $\Delta_v=2C(\pi)$ and $\Delta_v^2-4=4\bigl(C(\pi)^2-1\bigr)=4C'(\pi)S(\pi)$.
\end{proof}

\begin{lemma}[translation identity]\label{lem:shift}
For all $\lambda\in\C$,
\[
s_v(\pi,\lambda)=S(\pi,\lambda)\,C(\tau,\lambda)^2-C'(\pi,\lambda)\,S(\tau,\lambda)^2 .
\]
\end{lemma}

\begin{proof}
Let $Y=T^{-1}s_v$, that is, $Y=s_v(\cdot+\tau)$. Then $Y$ solves the Mathieu equation, with $Y(-\tau)=s_v(0)=0$ and $Y'(-\tau)=s_v'(0)=1$. The solution $C(-\tau)S-S(-\tau)C$ of the Mathieu equation vanishes at $-\tau$ and has derivative $W(C,S)=1$ there, so uniqueness gives $Y=C(-\tau)S-S(-\tau)C$. Inserting \eqref{eq:monodromy} at $x=-\tau$ into $s_v(\pi)=Y(\pi-\tau)=C(-\tau)S(\pi-\tau)-S(-\tau)C(\pi-\tau)$, we obtain
\[
s_v(\pi)=S(\pi)C(-\tau)^2+\bigl(S'(\pi)-C(\pi)\bigr)C(-\tau)S(-\tau)-C'(\pi)S(-\tau)^2.
\]
The middle term vanishes, since $C(\pi)=S'(\pi)$ by Lemma~\ref{lem:even}(a). Finally, $C(-\tau)^2=C(\tau)^2$ and $S(-\tau)^2=S(\tau)^2$ by \eqref{eq:parity}, and the identity of the lemma follows.
\end{proof}

Thus $\lambda$ is a Dirichlet eigenvalue of $L$ if and only if
\begin{equation}\label{eq:dirichlet}
S(\pi,\lambda)\,C(\tau,\lambda)^2=C'(\pi,\lambda)\,S(\tau,\lambda)^2.
\end{equation}
For the Mathieu operator on $[0,\pi]$, the zeros of $S(\pi,\cdot)$ and $C'(\pi,\cdot)$ are the Dirichlet and Neumann eigenvalues, with the odd solution $S$ and the even solution $C$ as eigenfunctions, respectively, by \eqref{eq:parity}. By Lemma~\ref{lem:even}(b), the union of these two sets is the union of the periodic and antiperiodic spectra. This splitting is also in \cite[\S4.5]{Vel25}.

\section{Estimates near \texorpdfstring{$n^2$}{n\^{}2}}\label{sec:asymptotics}

We now pass from the exact identities of Section~\ref{sec:shift} to estimates near $\lambda=n^2$, where $\mu_n$ and $\lambda_n^\pm$ lie for large $n$. In the rest of the paper, $Q>0$ is fixed, and $q\in\C$ with $|q|\le Q$ and $\tau\in[0,\pi)$ are arbitrary. Unless stated otherwise, $n$ is large, and all thresholds and constants, including the implied constants in $O$-terms, depend only on $Q$. Where a statement fixes $q$, or $a$ and $b$, the thresholds and constants may depend on them as well.

Set $U(x)=\int_0^xu(t)\,dt=q\sin2x$. For large $k$, the local frequency of the Mathieu equation with $\lambda=k^2$, that is, the square root of $k^2-u(t)$ near $k$, is $k-u(t)/(2k)+O(k^{-3})$, and integrating its first two terms from $0$ to $x$ gives the phase $\varphi=\varphi(x,k)=kx-U(x)/(2k)$. With this phase, the errors in the estimates below are $O(k^{-2})$, whereas $kx$ alone gives only $O(k^{-1})$.

\begin{lemma}[solution estimates]\label{lem:volterra}
For $x\in[0,\pi]$ and $k\in\C$ with $|k|\ge1$ and $|\Img k|\le1$,
\[
C(x,k^2)=\cos\varphi+O(k^{-2}),\qquad
kS(x,k^2)=\sin\varphi+O(k^{-2}),\qquad
-\frac{C'(x,k^2)}{k}=\sin\varphi+O(k^{-2}),
\]
uniformly in $x$ and $k$. Moreover, $|\cos\varphi|,|\sin\varphi|\le e^{\pi+Q/2}$.
\end{lemma}

\begin{proof}
Since $u=2q\cos2x$ and $U=q\sin2x$ with $|q|\le Q$, we have $|u|\le2Q$, $|u'|\le4Q$ and $|U|\le Q$ on $[0,\pi]$. Hence $|\Img\varphi|\le x|\Img k|+|U(x)|/(2|k|)\le\pi+Q/2$. The inequalities $|\cos z|,|\sin z|\le e^{|\Img z|}$ now give the stated bounds on $\cos\varphi$ and $\sin\varphi$. They also give $|\sin kt|,|\cos kt|\le e^{\pi}$ for $0\le t\le\pi$, as $|\Img(kt)|\le\pi$.

As $U'=u$, we have $\varphi'=k-u/(2k)$ and $\varphi''=-u'/(2k)$, and squaring gives $(\varphi')^2=k^2-u+u^2/(4k^2)$. Differentiating $\cos\varphi$ and $\sin\varphi$ twice therefore gives
\begin{align}
(\cos\varphi)''+(k^2-u)\cos\varphi
&=\frac{u'}{2k}\sin\varphi-\frac{u^2}{4k^2}\cos\varphi,
\label{eq:res-cos}\\
(\sin\varphi)''+(k^2-u)\sin\varphi
&=-\frac{u'}{2k}\cos\varphi-\frac{u^2}{4k^2}\sin\varphi.
\label{eq:res-sin}
\end{align}
Both right-hand sides have modulus at most $e^{\pi+Q/2}\bigl(2Q/|k|+Q^2/|k|^2\bigr)$, so they are $O(k^{-1})$, as $|k|\ge1$.

For $\eta=C-\cos\varphi$ or $\eta=kS-\sin\varphi$, let $h$ denote the right-hand side of \eqref{eq:res-cos} or \eqref{eq:res-sin}, respectively. Since $C$ and $kS$ solve the Mathieu equation $y''+(k^2-u)y=0$, subtraction gives $\eta''+(k^2-u)\eta=-h$, that is, $\eta''+k^2\eta=u\eta-h$. Moreover, $\varphi(0)=0$ because $U(0)=0$, so $\eta(0)=0$ in both cases. Differentiation gives $\eta'=C'+\varphi'\sin\varphi$ or $\eta'=kS'-\varphi'\cos\varphi$, so $\eta'(0)$ equals $0$ or $k-\varphi'(0)=u(0)/(2k)$, respectively. Variation of constants with these initial values gives
\begin{equation}
\eta(x)=\frac{\eta'(0)}{k}\sin kx
+\int_0^x\frac{\sin k(x-t)}{k}\bigl(u(t)\eta(t)-h(t)\bigr)\,dt.
\label{eq:eta-int}
\end{equation}
As $0\le x-t\le\pi$ in \eqref{eq:eta-int}, the bounds on $u$ and on $\sin kt$ yield
\[
|\eta(x)|\le\frac{e^{\pi}}{|k|}\Bigl(|\eta'(0)|+\int_0^x|h(t)|\,dt\Bigr)
+\frac{2Qe^{\pi}}{|k|}\int_0^x|\eta(t)|\,dt.
\]
Here $|\eta'(0)|$ and $\int_0^x|h(t)|\,dt$ are $O(k^{-1})$ uniformly for $0\le x\le\pi$, so the first term on the right-hand side is $O(k^{-2})$. As the coefficient $2Qe^{\pi}/|k|$ is at most $2Qe^{\pi}$, Gronwall's inequality gives $\eta=O(k^{-2})$ in both cases, and this proves the first two estimates.

For the third estimate, let $\eta=C-\cos\varphi$. Since $\eta'(0)=0$ and the kernel $\sin k(x-t)$ vanishes at $t=x$, differentiating \eqref{eq:eta-int} with respect to $x$ gives
\[
\eta'(x)=\int_0^x\cos k(x-t)\,\bigl(u(t)\eta(t)-h(t)\bigr)\,dt.
\]
Here $|\cos k(x-t)|\le e^{\pi}$, $u\eta=O(k^{-2})$ by the first estimate, and $h=O(k^{-1})$, so $\eta'=O(k^{-1})$. Since $u\sin\varphi$ is bounded in modulus by $2Qe^{\pi+Q/2}$, dividing $\eta'=C'+\bigl(k-u/(2k)\bigr)\sin\varphi$ by $-k$ gives
\[
-\frac{C'}{k}-\sin\varphi
=-\frac{\eta'}{k}-\frac{u}{2k^2}\sin\varphi
=O(k^{-2}).\qedhere
\]
\end{proof}

For $n\ge2$, put $D_n=\{\lambda:\ |\lambda-n^2|\le n\}$. This disc lies in the half-plane $\operatorname{Re}\lambda>0$, and we let $k=\sqrt\lambda$ be the principal square root there. On $D_n$, $|k+n|\ge n+\operatorname{Re}k\ge n$ gives $|k-n|=|\lambda-n^2|/|k+n|\le1$, and then $|k+n|=|2n+(k-n)|$ gives $3n/2\le2n-1\le|k+n|\le2n+1\le3n$. Hence
\begin{equation}\label{eq:kn}
|\lambda-n^2|/(3n)\le|k-n|\le2|\lambda-n^2|/(3n)\qquad\text{on }D_n.
\end{equation}
In particular, $|\Img k|\le|k-n|\le1$ and $|k|\ge n-1\ge n/2\ge1$, so Lemma~\ref{lem:volterra} applies to this $k$, and its errors $O(k^{-2})$ are $O(n^{-2})$ on $D_n$. Since $U(\pi)=0$, we have $\varphi(\pi,k)=k\pi$, so on $D_n$ the lemma at $x=\pi$ gives
\begin{equation}\label{eq:endpoint}
C(\pi)=\cos k\pi+O(n^{-2}),\qquad kS(\pi)=\sin k\pi+O(n^{-2}),\qquad -C'(\pi)/k=\sin k\pi+O(n^{-2}).
\end{equation}

The Dirichlet and periodic or antiperiodic eigenvalues near $n^2$ are counted in \cite{DM06,DM20}. We locate them within $O(n^{-1})$ of $n^2$, uniformly in $q$ and $\tau$, and identify $\lambda_n^\pm$ as the zeros in $D_n$ of $S(\pi,\cdot)$ and $C'(\pi,\cdot)$.

\begin{lemma}[localization]\label{lem:local}
There are $n_0$ and $M$ such that for $n\ge n_0$, each of $S(\pi,\cdot)$, $C'(\pi,\cdot)$ and $s_v(\pi,\cdot)$ has exactly one zero $\lambda$ in $D_n$. This zero is simple and satisfies $|\lambda-n^2|\le M/n\le1$ and $|\sqrt\lambda-n|\le M/n^2$.
\end{lemma}

\begin{proof}
On $D_n$, Lemma~\ref{lem:volterra} at $x=\tau$ gives $C(\tau)=\cos\varphi(\tau,k)+O(n^{-2})$ and $kS(\tau)=\sin\varphi(\tau,k)+O(n^{-2})$, and it bounds $|\cos\varphi(\tau,k)|$, $|\sin\varphi(\tau,k)|$ and $|\sin k\pi|=|\sin\varphi(\pi,k)|$ by $e^{\pi+Q/2}$. Hence products of these estimates have errors $O(n^{-2})$. Multiplying the identity of Lemma~\ref{lem:shift} by $k$ and inserting these estimates and \eqref{eq:endpoint} gives
\begin{align*}
k\,s_v(\pi)&=kS(\pi)\,C(\tau)^2-\frac{C'(\pi)}{k}\,\bigl(kS(\tau)\bigr)^2\\
&=\sin k\pi\,\bigl(\cos^2\varphi(\tau,k)+\sin^2\varphi(\tau,k)\bigr)+O(n^{-2})=\sin k\pi+O(n^{-2}).
\end{align*}
Thus there is a constant $K_1$ such that, if $F$ denotes any of the functions $kS(\pi)$, $-C'(\pi)/k$ and $k\,s_v(\pi)$, then
\begin{equation}\label{eq:rouche}
|F-\sin k\pi|\le K_1n^{-2}\qquad\text{on }D_n.
\end{equation}

To compare $F$ with $\sin k\pi$ by Rouch\'e's theorem, write $k=n+w$, so that $\sin k\pi=(-1)^n\sin\pi w$ and $\lambda-n^2=w(k+n)$. By \eqref{eq:kn}, $|w|\le2/3$ on $D_n$ and $|w|\ge1/3$ on $\partial D_n$. The function $\sin(\pi w)/w$, defined to be $\pi$ at $w=0$, is continuous and nonzero on $|w|\le2/3$, as $\sin\pi w$ vanishes only at the integers. Hence there is an absolute constant $c_0>0$ such that $|\sin k\pi|=|\sin\pi w|\ge c_0|w|$ on $D_n$, so $|\sin k\pi|\ge c_0/3$ on $\partial D_n$. Moreover,
\[
\sin k\pi=(\lambda-n^2)\,\frac{(-1)^n}{k+n}\,\frac{\sin\pi w}{w}.
\]
The factor multiplying $\lambda-n^2$ is holomorphic and nonzero throughout $D_n$, so $\lambda=n^2$ is the only zero of $\sin k\pi$ in $D_n$, and it is simple.

As $k$ is holomorphic and nonzero near $D_n$, and $S(\pi,\cdot)$, $C'(\pi,\cdot)$ and $s_v(\pi,\cdot)$ are entire by Section~\ref{sec:shift}, $F$ is holomorphic near $D_n$, with the same zeros and the same orders as the corresponding one of these three functions. If $K_1n^{-2}<c_0/3$, then \eqref{eq:rouche} and the lower bound for $\sin k\pi$ on $\partial D_n$ give
\[
|F-\sin k\pi|\le K_1n^{-2}<c_0/3\le|\sin k\pi|\qquad\text{on }\partial D_n.
\]
Hence $F$ has no zeros on $\partial D_n$, and by Rouch\'e's theorem $F$ and $\sin k\pi$ have the same number of zeros in $D_n$, counted with multiplicity, namely one. Therefore each of $S(\pi,\cdot)$, $C'(\pi,\cdot)$ and $s_v(\pi,\cdot)$ has exactly one zero in $D_n$, and this zero is simple.

It remains to locate this zero $\lambda$. As the corresponding $F$ vanishes there, the lower bound for $\sin k\pi$ and \eqref{eq:rouche} give $c_0|w|\le|\sin k\pi|\le K_1n^{-2}$, so $|\sqrt\lambda-n|=|w|\le K_1/(c_0n^2)$. By \eqref{eq:kn}, $|\lambda-n^2|\le3n|w|\le3K_1/(c_0n)$. Thus the lemma holds with $M=3K_1/c_0$ and any $n_0\ge\max\{2,M\}$ with $K_1n_0^{-2}<c_0/3$.
\end{proof}

Denote by $\ls$, $\lc$ and $\mu_n$ the respective zeros in $D_n$ of $S(\pi,\cdot)$, $C'(\pi,\cdot)$ and $s_v(\pi,\cdot)$ given by Lemma~\ref{lem:local}. By Section~\ref{sec:shift}, $\mu_n$ is the only Dirichlet eigenvalue of $L$ in $D_n$.

\begin{lemma}[identification of the edges]\label{lem:edges}
There is $n_0$ such that for $n\ge n_0$, the periodic ($n$~even) or antiperiodic ($n$~odd) eigenvalues of $L$ in $D_n$, counted with algebraic multiplicity, are $\ls$ and $\lc$. Moreover, for $x\in\R$,
\[
S(x+\pi,\ls)=(-1)^nS(x,\ls),\qquad C(x+\pi,\lc)=(-1)^nC(x,\lc).
\]
\end{lemma}

\begin{proof}
Let $\lambda$ be $\ls$ or $\lc$. Then $C'(\pi,\lambda)S(\pi,\lambda)=0$, so $C(\pi,\lambda)^2=1$ by Lemma~\ref{lem:even}(a). Also, $|k-n|\le M/n^2$ by Lemma~\ref{lem:local}, so $\cos k\pi=(-1)^n\cos\pi(k-n)=(-1)^n+O(n^{-4})$, and \eqref{eq:endpoint} gives $C(\pi,\lambda)=(-1)^n+O(n^{-2})$. As $C(\pi,\lambda)=\pm1$, this forces $C(\pi,\lambda)=(-1)^n$ for large $n$, and then $S'(\pi,\lambda)=(-1)^n$ by Lemma~\ref{lem:even}(a). At $\ls$, where $S(\pi)=0$, the second identity of \eqref{eq:monodromy} therefore reduces to $S(x+\pi)=(-1)^nS(x)$, and at $\lc$, where $C'(\pi)=0$, the first reduces to $C(x+\pi)=(-1)^nC(x)$.

For the first assertion, Lemma~\ref{lem:even}(b) gives $\Delta_v=2C(\pi)=2(-1)^n$ at $\ls$ and $\lc$, so by Section~\ref{sec:shift} both are periodic eigenvalues of $L$ if $n$ is even and antiperiodic eigenvalues if $n$ is odd. Moreover, $\Delta_v^2-4=4C'(\pi)S(\pi)$ has no other zeros in $D_n$ by Lemma~\ref{lem:local}, so $L$ has no other periodic or antiperiodic eigenvalues there. It remains to determine the multiplicities. With the $L^2([0,\pi])$ norm of \cite{DM06}, normalized so that $\|e^{2ix}\|=1$, the potential $v$ has norm $(|a|^2+|b|^2)^{1/2}=\sqrt2\,|q|$. By \cite[\S2.2]{DM06}, for $n>6\sqrt2\,|q|$ the periodic ($n$~even) or antiperiodic ($n$~odd) problem for $L$ has exactly two eigenvalues, counted with algebraic multiplicity, in the disc $|\lambda-n^2|\le n/2$. This disc lies in $D_n$ and contains $\ls$ and $\lc$ by Lemma~\ref{lem:local}. Hence $\ls$ and $\lc$ each have algebraic multiplicity one if they are distinct, and their common value has algebraic multiplicity two if they coincide.
\end{proof}

To bring the edge deviations $\mu_n-\ls$ and $\mu_n-\lc$ into \eqref{eq:dirichlet}, we write
\begin{equation}\label{eq:factor}
S(\pi,\lambda)=(\lambda-\ls)\,f_n(\lambda),\qquad C'(\pi,\lambda)=(\lambda-\lc)\,g_n(\lambda).
\end{equation}
By Lemma~\ref{lem:local}, $\ls$ and $\lc$ are the only zeros in $D_n$ of $S(\pi,\cdot)$ and $C'(\pi,\cdot)$, respectively, and both are simple. Hence $f_n$ and $g_n$ are entire and have no zeros in $D_n$. Put $r_n=-g_n/f_n$ on $D_n$.

\begin{lemma}[ratio of the factors]\label{lem:factor}
As $n\to\infty$, uniformly for $|\lambda-n^2|\le1$, we have $r_n(\lambda)=n^2\bigl(1+O(n^{-2})\bigr)$.
\end{lemma}

\begin{proof}
We compare $f_n$ and $-g_n/\lambda$ with their common value for $q=0$. For $q=0$ we have $S=\sin(kx)/k$ and $C'=-k\sin kx$, so $S(\pi,\lambda)$ and $-C'(\pi,\lambda)/\lambda$ both equal $S^0(\lambda)=\sin(k\pi)/k$. The function $S^0$ is entire, since $\sin(k\pi)/k$ is a power series in $k^2=\lambda$ that converges for all $\lambda$, and $S^0(n^2)=0$. Hence $f_n^0=S^0/(\lambda-n^2)$ is entire as well. Dividing the last two estimates of \eqref{eq:endpoint} by $k$, with $|k|\ge n/2$, gives
\[
S(\pi,\lambda)-S^0(\lambda)=O(n^{-3}),\qquad-\frac{C'(\pi,\lambda)}{\lambda}-S^0(\lambda)=O(n^{-3})\qquad\text{on }D_n.
\]
By the factorization \eqref{eq:factor} and $S^0=(\lambda-n^2)f_n^0$,
\begin{align*}
(\lambda-\ls)\bigl(f_n-f_n^0\bigr)&=S(\pi,\lambda)-S^0(\lambda)+(\ls-n^2)f_n^0,\\
(\lambda-\lc)\Bigl(-\frac{g_n}{\lambda}-f_n^0\Bigr)&=-\frac{C'(\pi,\lambda)}{\lambda}-S^0(\lambda)+(\lc-n^2)f_n^0.
\end{align*}

On $\partial D_n$ we have $|\lambda-n^2|=n$ and $|S^0|\le2e^\pi/n$, so $f_n^0=O(n^{-2})$. Moreover, $|\ls-n^2|\le M/n$ and $|\lc-n^2|\le M/n$ by Lemma~\ref{lem:local}. Hence both right-hand sides are $O(n^{-3})$ on $\partial D_n$, while $|\lambda-\ls|$ and $|\lambda-\lc|$ are at least $n-M/n\ge n/2$ there. Thus $f_n-f_n^0$ and $-g_n/\lambda-f_n^0$ are $O(n^{-4})$ on $\partial D_n$. As $0\notin D_n$, these two functions are holomorphic near $D_n$, and the maximum modulus principle gives the same bounds on $D_n$.

Finally, let $|\lambda-n^2|\le1$. Then $w=k-n=(\lambda-n^2)/(k+n)$ satisfies $|w|\le1/n$ by \eqref{eq:kn}, so $k=n\bigl(1+O(n^{-2})\bigr)$. Since $\sin k\pi=(-1)^n\pi w\bigl(1+O(w^2)\bigr)$ and $\lambda-n^2=w(k+n)$,
\[
f_n^0(\lambda)=\frac{(-1)^n\pi\bigl(1+O(n^{-2})\bigr)}{k(k+n)}=(-1)^n\frac{\pi}{2n^2}\bigl(1+O(n^{-2})\bigr).
\]
For large $n$ this gives $|f_n^0|\ge\pi/(4n^2)$, so the two bounds $O(n^{-4})$ on $D_n$ become $f_n=f_n^0\bigl(1+O(n^{-2})\bigr)$ and $-g_n/\lambda=f_n^0\bigl(1+O(n^{-2})\bigr)$. Dividing the second by the first gives $r_n/\lambda=1+O(n^{-2})$, and the estimate follows, as $\lambda=n^2\bigl(1+O(n^{-2})\bigr)$.
\end{proof}

\section{Asymptotics of the signed Mathieu gap}\label{sec:gap}

Recall from Section~\ref{sec:shift} that $L$ and the Mathieu operator $-d^2/dx^2+u$, where $u(x)=2q\cos2x$ and $q^2=ab$, have the same spectral gaps. By Lemma~\ref{lem:edges}, the edges $\lambda_n^\pm$ are $\ls$ and $\lc$ in some order. Unlike in \cite{AD12,DM20}, where $\gamma_n$ denotes $\lambda_n^+-\lambda_n^-$, we put $\gamma_n=\lc-\ls$. Then $\gamma_n=\pm(\lambda_n^+-\lambda_n^-)$, so $|\gamma_n|$ is the gap length, and the sign records which edge belongs to $C$ and which to $S$. We need this sign, since $S(\pi,\cdot)$ and $C'(\pi,\cdot)$ enter \eqref{eq:dirichlet} with different coefficients. Moreover, $\lambda_n^-+\lambda_n^+=\ls+\lc$, so the midpoint of the two edges is $\ls+\frac12\gamma_n$ and
\begin{equation}\label{eq:delta-mid}
\delta_n=(\mu_n-\ls)-\tfrac12\gamma_n.
\end{equation}

For the Mathieu operator with real $q$, Harrell \cite{Har81} found the leading term of the gap length up to a constant factor, and Avron and Simon \cite{AS81} found this term to be $8|q|^nF_n$. For fixed complex $q\neq0$, Anahtarci and Djakov \cite{AD12} proved, up to sign, a refined expansion of $\lambda_n^+-\lambda_n^-$, sharpened in~\cite{Ana26}. The sign is open in both, because their reduction to $\operatorname{span}\{e^{-inx},e^{inx}\}$ splits the problem into two scalar equations without deciding which edge solves which. The parity of $C$ and $S$ decides it, giving the signed expansion of $\gamma_n$ to the order we need.

\begin{proposition}[signed gap asymptotics]\label{lem:gap}
Let $q\in\C\setminus\{0\}$ be fixed. As $n\to\infty$,
\begin{equation}\label{eq:gap-lead}
\gamma_n=8q^nF_n\Bigl(1-\frac{q^2}{4n^3}+O(n^{-4})\Bigr).
\end{equation}
\end{proposition}

\begin{proof}
Let $y=C(\cdot,\lc)$ and $z=z_c=\lc-n^2$, or $y=S(\cdot,\ls)$ and $z=z_s=\ls-n^2$. By Lemma~\ref{lem:edges}, $y(x+\pi)=(-1)^ny(x)$, so $y$ is an eigenfunction of the Mathieu operator with eigenvalue $n^2+z$, periodic for even~$n$ and antiperiodic for odd~$n$. We work in $L^2([0,\pi])$ with the orthogonal basis $e^{ijx}$, $j\in n+2\Z$, write $y_j$ for the coefficients of $y$, and let $\Pi$ be the orthogonal projection onto $\operatorname{span}\{e^{-inx},e^{inx}\}$. By Lemma~\ref{lem:local}, $|z|\le M/n$, so $|z|<1$ for large $n$. Let $G_z$ act as $0$ on $e^{\pm inx}$ and as $(n^2-j^2+z)^{-1}$ on $e^{ijx}$ for ${j\neq\pm n}$. Since $|n^2-j^2|\ge4n-4$ for these $j$, $\|G_z\|\le1/(4n-5)$, and multiplication by $u$ has norm at most $2|q|$, so $\|G_zu\|\le2|q|/(4n-5)<1$. Put $p=\Pi y$. In coordinates, the eigenvalue equation $y''+(n^2+z)y=uy$ reads $(n^2-j^2+z)y_j=(uy)_j$. For~${j\neq\pm n}$ this says $y-p=G_zuy$, so $y=(I-G_zu)^{-1}p$, and $p\neq0$, as $y\neq0$. For~${j=\pm n}$ it says $zp=\Pi uy$, and hence $zp=\Pi u(I-G_zu)^{-1}p$, as in the reduction of \cite[\S\S2.1--2.2]{DM06}.

For $j\in n+2\Z$ the product $y(x)e^{-ijx}$ is $\pi$-periodic, so $y_j=\frac1\pi\int_{-\pi/2}^{\pi/2}y(x)e^{-ijx}\,dx$. As $C$ is even and $S$ is odd by \eqref{eq:parity}, the substitution $x\mapsto-x$ gives $y_{-n}=y_n$ for $y=C(\cdot,\lc)$ and $y_{-n}=-y_n$ for $y=S(\cdot,\ls)$. Thus $p$ is a nonzero multiple of $\cos nx$ or of $\sin nx$, respectively. By \cite[\S2]{AD12}, where the parameter of the Mathieu potential is our $q$, the matrix of $\Pi u(I-G_zu)^{-1}\Pi$ in the basis $e^{-inx},e^{inx}$ has the form $\bigl(\begin{smallmatrix}\alpha_n(z)&\beta_n(z)\\ \beta_n(z)&\alpha_n(z)\end{smallmatrix}\bigr)$, so this operator multiplies $\cos nx$ by $\alpha_n(z)+\beta_n(z)$ and $\sin nx$ by $\alpha_n(z)-\beta_n(z)$. Hence $zp=\Pi u(I-G_zu)^{-1}p$ gives
\[
z_c=\alpha_n(z_c)+\beta_n(z_c),\qquad z_s=\alpha_n(z_s)-\beta_n(z_s).
\]
Subtracting the second equation from the first, with $z_c-z_s=\lc-\ls=\gamma_n$, gives
\[
\gamma_n-\bigl(\alpha_n(z_c)-\alpha_n(z_s)\bigr)=\beta_n(z_c)+\beta_n(z_s).
\]
Since $\{z_c,z_s\}=\{\lambda_n^--n^2,\lambda_n^+-n^2\}$, the estimates of \cite[\S\S3--4]{AD12} give
\[
\alpha_n(z_c)-\alpha_n(z_s)=-\gamma_n\Bigl(\frac{q^2}{8n^2}+O(n^{-4})\Bigr),\qquad\beta_n(z)=4q^nF_n\Bigl(1+\frac{q^2}{8n^2}-\frac{q^2}{4n^3}+O(n^{-4})\Bigr)\quad(z=z_c,z_s).
\]
Inserting these estimates into the identity for $\gamma_n$ gives
\[
\gamma_n\Bigl(1+\frac{q^2}{8n^2}+O(n^{-4})\Bigr)=8q^nF_n\Bigl(1+\frac{q^2}{8n^2}-\frac{q^2}{4n^3}+O(n^{-4})\Bigr).
\]
Dividing by the bracket on the left removes the term $q^2/(8n^2)$ on the right up to $O(n^{-4})$ and gives \eqref{eq:gap-lead}.
\end{proof}

\section{Position of the Dirichlet eigenvalue and the main expansion}\label{sec:main}

We now locate $\mu_n$ relative to the edges $\ls$ and $\lc$, and deduce the expansion of $\delta_n$. By \eqref{eq:dirichlet}, the position of $\mu_n$ is governed by the coefficients $C(\tau,\mu_n)^2$ and $S(\tau,\mu_n)^2$. Lemma~\ref{lem:volterra} describes these through the phase at $x=\tau$. Since $\sqrt{\mu_n}=n+O(n^{-2})$ by Lemma~\ref{lem:local}, replacing $\sqrt{\mu_n}$ by $n$ changes this phase by $O(n^{-2})$, and we put
\begin{equation}\label{eq:theta}
\Th=\varphi(\tau,n)=n\tau-\frac{q\sin2\tau}{2n}=n\tau+\frac{i(a-b)}{4n},\qquad |\Img\Th|\le\frac{|q|}{2n}.
\end{equation}
The last equality follows from \eqref{eq:shift-cons}, and the bound holds because $\tau$ is real.

\needspace{8\baselineskip}
\begin{proposition}[position estimate]\label{prop:position}
For all large $n$, there is $w_n\in\C$ with
\[
\mu_n=(1-w_n)\ls+w_n\lc,\qquad\delta_n=\gamma_n\bigl(w_n-\tfrac12\bigr).
\]
As $n\to\infty$, uniformly in $\tau$ and in $q$ on bounded subsets of $\C$,
\begin{equation}\label{eq:weight}
w_n=\bigl(\sin\Th+O(n^{-2})\bigr)^2\bigl(1+O(n^{-2})\bigr),\qquad 1-w_n=\bigl(\cos\Th+O(n^{-2})\bigr)^2\bigl(1+O(n^{-2})\bigr).
\end{equation}
In particular, $w_n=\sin^2\Th+O(n^{-2})$ and
\begin{equation}\label{eq:position}
\mu_n=\ls\cos^2\Th+\lc\sin^2\Th+O(\gamma_nn^{-2}),\qquad \delta_n=-\tfrac12\gamma_n\cos2\Th+O(\gamma_nn^{-2}).
\end{equation}
\end{proposition}

\begin{proof}
Let $\mu=\mu_n$ and $k=\sqrt\mu$. As $f_n$ has no zeros in $D_n$, inserting the factorization \eqref{eq:factor} into the Dirichlet condition \eqref{eq:dirichlet} at $\lambda=\mu$ and dividing by $f_n(\mu)$, with $r_n=-g_n/f_n$, gives
\begin{equation}\label{eq:exact}
(\mu-\ls)\,C(\tau,\mu)^2+(\mu-\lc)\,r_n(\mu)\,S(\tau,\mu)^2=0.
\end{equation}
To solve \eqref{eq:exact} for $\mu-\ls$, we estimate its two coefficients and show that their sum is $1+O(n^{-2})$.

By Lemma~\ref{lem:local}, $k-n=O(n^{-2})$ and $|\mu-n^2|\le1$. The phase of Lemma~\ref{lem:volterra} at $x=\tau$ and $\lambda=\mu$ is $\varphi(\tau,k)$, and it differs from $\Th=\varphi(\tau,n)$ by
\[
d=\varphi(\tau,k)-\Th=(k-n)\tau-\frac{U(\tau)}{2}\Bigl(\frac1k-\frac1n\Bigr)=(k-n)\Bigl(\tau+\frac{U(\tau)}{2kn}\Bigr)=O(n^{-2}),
\]
as $0\le\tau<\pi$, $|U(\tau)|\le Q$ and $|kn|\ge n^2/2$. Hence, by Lemma~\ref{lem:volterra} with the phase $\varphi(\tau,k)=\Th+d$,
\begin{align*}
C(\tau,\mu)&=\cos(\Th+d)+O(n^{-2})=\cos\Th+O(n^{-2}),\\
kS(\tau,\mu)&=\sin(\Th+d)+O(n^{-2})=\sin\Th+O(n^{-2}).
\end{align*}
where the second equalities follow from the addition formulas, as $\cos d=1+O(n^{-4})$, $\sin d=O(n^{-2})$ and $|\sin\Th|,|\cos\Th|\le e^{Q}$ by \eqref{eq:theta}. Since $r_n(\mu)=n^2\bigl(1+O(n^{-2})\bigr)$ by Lemma~\ref{lem:factor} and $k=n+O(n^{-2})$, we have $r_n(\mu)/k^2=1+O(n^{-2})$. Squaring the estimates for $C(\tau,\mu)$ and $kS(\tau,\mu)$ therefore gives
\begin{equation}\label{eq:sq}
\begin{aligned}
C(\tau,\mu)^2&=\bigl(\cos\Th+O(n^{-2})\bigr)^2,\\
r_n(\mu)\,S(\tau,\mu)^2&=\frac{r_n(\mu)}{k^2}\bigl(kS(\tau,\mu)\bigr)^2=\bigl(\sin\Th+O(n^{-2})\bigr)^2\bigl(1+O(n^{-2})\bigr).
\end{aligned}
\end{equation}
Expanding the squares in \eqref{eq:sq} and adding gives
\begin{equation}\label{eq:sq-exp}
\sigma=C(\tau,\mu)^2+r_n(\mu)\,S(\tau,\mu)^2=\cos^2\Th+\sin^2\Th+O(n^{-2})=1+O(n^{-2}),
\end{equation}
so $\sigma\neq0$ for large $n$. Since $\mu-\lc=\mu-\ls-\gamma_n$, \eqref{eq:exact} becomes $(\mu-\ls)\,\sigma=\gamma_n\,r_n(\mu)\,S(\tau,\mu)^2$. Dividing by $\sigma$ gives $\mu-\ls=w_n\gamma_n$, that is, $\mu=(1-w_n)\ls+w_n\lc$, where
\begin{equation}\label{eq:quot}
w_n=\frac{r_n(\mu)\,S(\tau,\mu)^2}{\sigma},\qquad 1-w_n=\frac{C(\tau,\mu)^2}{\sigma},
\end{equation}
and then \eqref{eq:delta-mid} gives $\delta_n=\gamma_n\bigl(w_n-\frac12\bigr)$. By \eqref{eq:sq-exp}, $1/\sigma=1+O(n^{-2})$, so \eqref{eq:sq} and \eqref{eq:quot} give \eqref{eq:weight}. As $\sin\Th$ is bounded, \eqref{eq:weight} gives $w_n=\sin^2\Th+O(n^{-2})$. Inserting this estimate into $\mu_n=\ls+w_n\gamma_n$ and $\delta_n=\gamma_n\bigl(w_n-\frac12\bigr)$ gives \eqref{eq:position}, since $\ls+\gamma_n\sin^2\Th=\ls\cos^2\Th+\lc\sin^2\Th$ and $\sin^2\Th-\frac12=-\frac12\cos2\Th$.
\end{proof}

If $\gamma_n\neq0$, then $|\mu_n-\ls|/|\gamma_n|=|w_n|$ and $|\mu_n-\lc|/|\gamma_n|=|1-w_n|$, and by Proposition~\ref{prop:position} these are $|\sin\Th|^2+O(n^{-2})$ and $|\cos\Th|^2+O(n^{-2})$. They are the two ratios $|\mu_n-\lambda_n^\pm|/|\lambda_n^+-\lambda_n^-|$ in the Riesz-basis criteria of \cite{DM12,GT12}, and by \eqref{eq:theta} they are bounded, in agreement with \cite{DM11}.

In the self-adjoint case $b=\bar a$, the parameter $q$ is real, as $q^2=|a|^2$, and so are the eigenvalue $\mu_n$ and the edges $\ls$ and $\lc$. Then $C(\tau,\mu_n)$ and $S(\tau,\mu_n)$ are real, $r_n(\mu_n)$ is real by \eqref{eq:factor}, and $r_n(\mu_n)>0$ by Lemma~\ref{lem:factor}. Both coefficients in \eqref{eq:exact} are therefore nonnegative, so \eqref{eq:quot} gives $0\le w_n\le1$, and $\mu_n$ lies in the closed gap, at relative position $w_n=\sin^2\Th+O(n^{-2})$ from $\ls$ toward $\lc$ when $\gamma_n\neq0$.

\begin{remark}
The change $(q,\tau)\mapsto(-q,\tau+\pi/2)$ of Section~\ref{sec:shift} replaces $u$ by $u(\cdot+\pi/2)$, and by \eqref{eq:theta} it replaces $\Th$ by $\Th+n\pi/2$ modulo $n\pi$. By Lemma~\ref{lem:edges} and \eqref{eq:parity}, translation by $\pi/2$ turns $S(\cdot,\ls)$ and $C(\cdot,\lc)$ into eigenfunctions for the translated potential with the same parity for even $n$ and the opposite parity for odd $n$. For even $n$ nothing in \eqref{eq:position} changes. For odd $n$, $\ls$ and $\lc$ are interchanged and $\gamma_n$ changes sign, while $\sin^2\Th$ and $\cos^2\Th$ are interchanged and $\cos2\Th$ changes sign, so \eqref{eq:position} is invariant.
\end{remark}

\needspace{4\baselineskip}
Combining \eqref{eq:position} with Proposition~\ref{lem:gap} gives the expansion of $\delta_n$ in terms of $a$ and $b$.

\begin{theorem}[main expansion]\label{thm:main}
Let $q\in\C\setminus\{0\}$. As $n\to\infty$, uniformly in $a,b\in\C$ with $ab=q^2$ and $|a|=|b|$,
\begin{equation}\label{eq:main2}
\delta_n=-2F_n\,(a^n+b^n)+\frac{F_n}{n}\,(a-b)(a^n-b^n)+O\Bigl(\frac{F_n|a|^n}{n^2}\Bigr).
\end{equation}
\end{theorem}

\begin{proof}
By \eqref{eq:gap-lead}, $\gamma_n=8q^nF_n+O(F_n|q|^nn^{-3})$, and by \eqref{eq:theta}, $|\cos2\Th|\le e^{|\Img2\Th|}\le e^{|q|}$. Substituting the gap estimate into the second formula in \eqref{eq:position} therefore gives, uniformly in $\tau$,
\begin{equation}\label{eq:delta-q}
\delta_n=-4q^nF_n\cos2\Th+O\bigl(F_n|q|^nn^{-2}\bigr).
\end{equation}

We now express $q^n\cos2\Th$ through $a$ and $b$. By \eqref{eq:theta} we have $2i(\Th-n\tau)=-(a-b)/(2n)$. Writing $2\cos2\Th=e^{2i\Th}+e^{-2i\Th}$ and using $q^ne^{2in\tau}=a^n$ and $q^ne^{-2in\tau}=b^n$ from \eqref{eq:shift-cons}, we obtain
\begin{equation}\label{eq:exp-id}
\xi_n=\frac{a-b}{2n},\qquad 2q^n\cos2\Th=a^ne^{-\xi_n}+b^ne^{\xi_n},\qquad |\xi_n|\le\frac{|q|}{n}.
\end{equation}
The bound on $\xi_n$ holds because $|a-b|\le|a|+|b|=2|q|$. It gives $|e^{\pm\xi_n}-1\mp\xi_n|\le|\xi_n|^2e^{|\xi_n|}\le|q|^2e^{|q|}n^{-2}$, and with $|a|=|b|=|q|$ the identity in \eqref{eq:exp-id} becomes
\[
2q^n\cos2\Th=a^n+b^n-\frac{(a-b)(a^n-b^n)}{2n}+O\bigl(|a|^nn^{-2}\bigr).
\]
Multiplying by $-2F_n$ and inserting the result into \eqref{eq:delta-q}, we arrive at \eqref{eq:main2}. All estimates are uniform in $\tau$, and by \eqref{eq:shift} the pairs $(a,b)$ of the statement correspond to $\tau\in[0,\pi)$, so the expansion is uniform in $a$ and $b$.
\end{proof}

Along a sequence of indices with $|(a/b)^n-1|=O(n^{-1})$, the second term of \eqref{eq:main2} is $O(F_n|a|^nn^{-2})$, the size of the bound for the remainder. On such a sequence the expansion can be refined by one term.

\begin{theorem}[three-term expansion]\label{thm:refined}
Let $a,b\in\C\setminus\{0\}$ be fixed with $|a|=|b|$. As $n\to\infty$ along any sequence of indices with $|(a/b)^n-1|=O(n^{-1})$,
\begin{equation}\label{eq:main3}
\delta_n=-2F_n(a^n+b^n)+\frac{F_n}{n}(a-b)(a^n-b^n)-\frac{(a-b)^2a^nF_n}{2n^2}+O\bigl(F_n|a|^nn^{-3}\bigr),
\end{equation}
with an implied constant that may depend on the sequence.
\end{theorem}

\begin{proof}
By \eqref{eq:shift-cons}, $(a/b)^n-1=e^{4in\tau}-1=2ie^{2in\tau}\sin2n\tau$, and $\tau$ is real, so $|\sin2n\tau|=\frac12|(a/b)^n-1|=O(n^{-1})$ along the sequence. The correction $2\Th-2n\tau$ is $O(n^{-1})$ by \eqref{eq:theta}, so it has sine $O(n^{-1})$ and bounded cosine. As $|\cos2n\tau|\le1$, the addition formula gives $\sin2\Th=O(n^{-1})$.

We now show that $w_n=\sin^2\Th+O(n^{-3})$ along the sequence. Since $\sin\Th$ and $\cos\Th$ are bounded by \eqref{eq:theta}, expanding the squares in the two formulas of \eqref{eq:weight} gives $w_n-\sin^2\Th=O(n^{-2}|\sin\Th|+n^{-4})$ and $1-w_n-\cos^2\Th=O(n^{-2}|\cos\Th|+n^{-4})$. The two left-hand sides differ only in sign, so
\[
w_n-\sin^2\Th=O\bigl(n^{-2}\min\{|\sin\Th|,|\cos\Th|\}+n^{-4}\bigr).
\]
Since $|\sin\Th|^2+|\cos\Th|^2\ge|\sin^2\Th+\cos^2\Th|=1$, the larger of $|\sin\Th|$ and $|\cos\Th|$ is at least $1/\sqrt2$. Hence the smaller is at most $\sqrt2\,|\sin\Th\cos\Th|=|\sin2\Th|/\sqrt2=O(n^{-1})$, and the estimate follows.

Thus $w_n-\frac12=-\frac12\cos2\Th+O(n^{-3})$ along the sequence. Inserting this and \eqref{eq:gap-lead} into $\delta_n=\gamma_n\bigl(w_n-\frac12\bigr)$ from Proposition~\ref{prop:position}, with $|\cos2\Th|\le e^{|q|}$ by \eqref{eq:theta}, we obtain
\begin{equation}\label{eq:delta-q3}
\delta_n=-4q^nF_n\cos2\Th+O\bigl(F_n|q|^nn^{-3}\bigr).
\end{equation}

It remains to expand $q^n\cos2\Th$ to this order. By \eqref{eq:exp-id}, $|\xi_n|\le|q|/n$, so $e^{\pm\xi_n}=1\pm\xi_n+\frac12\xi_n^2+O(n^{-3})$, and inserting this into the identity in \eqref{eq:exp-id} gives
\begin{equation}\label{eq:exp-second}
2q^n\cos2\Th=a^n+b^n-\xi_n(a^n-b^n)+\tfrac12\xi_n^2(a^n+b^n)+O\bigl(|a|^nn^{-3}\bigr).
\end{equation}
Along the sequence, $b^n-a^n=b^n\bigl(1-(a/b)^n\bigr)=O(|a|^nn^{-1})$. Thus, by \eqref{eq:exp-id}, the quadratic term in \eqref{eq:exp-second} is
\[
\frac12\xi_n^2(a^n+b^n)=\frac{(a-b)^2}{8n^2}\bigl(2a^n+(b^n-a^n)\bigr)=\frac{(a-b)^2a^n}{4n^2}+O\bigl(|a|^nn^{-3}\bigr).
\]
Inserting this into \eqref{eq:exp-second}, multiplying by $-2F_n$ and using \eqref{eq:delta-q3}, we obtain \eqref{eq:main3}, as $2F_n\xi_n=F_n(a-b)/n$.
\end{proof}

Let $P_n=-2F_n(a^n+b^n)+\frac{F_n}{n}(a-b)(a^n-b^n)$ be the sum of the first two terms in \eqref{eq:main2}, so that $\delta_n-P_n$ is the remainder of this expansion. Theorems~\ref{thm:main} and~\ref{thm:refined}, and for $a=b$ Proposition~\ref{lem:gap}, give the following.

\begin{corollary}[sharpness of the main expansion]\label{cor:sharp}
Let $a,b\in\C\setminus\{0\}$ be fixed with $|a|=|b|$.

\noindent\textup{(a)} If $a\neq b$, then the bounds
\[
\bigl|\delta_n+2F_n(a^n+b^n)\bigr|\ge\frac{|a-b|^2F_n|a|^{n-1}}{4n},\qquad
|\delta_n-P_n|\ge\frac{|a-b|^2F_n|a|^n}{4n^2}
\]
hold, each for infinitely many $n$.

\noindent\textup{(b)} If $a=b$, then, as $n\to\infty$,
\[
\delta_n-P_n=\frac{a^{n+2}F_n}{n^3}+O\bigl(F_n|a|^nn^{-4}\bigr).
\]

\noindent Thus $\delta_n-P_n=o(F_n|a|^nn^{-2})$ as $n\to\infty$ if and only if $a=b$.
\end{corollary}

\begin{proof}
(a) Put $z=a/b$, so that $|z|=1$, $|a-b|=|a|\,|z-1|$ and $|a^n-b^n|=|a|^n|z^n-1|$. For every $n$,
\[
|z^{n+1}-1|+|z^n-1|\ge|z^{n+1}-z^n|=|z-1|,
\]
so for every $n$ at least one of $|z^n-1|$ and $|z^{n+1}-1|$ is at least $\frac12|z-1|$. Hence $|z^n-1|\ge\frac12|z-1|$ for infinitely many $n$. For these $n$, the second term of \eqref{eq:main2} has modulus
\[
\frac{F_n}{n}\,|a-b|\,|a|^n|z^n-1|\ge\frac{|a-b|^2F_n|a|^{n-1}}{2n}.
\]
By \eqref{eq:main2} and the reverse triangle inequality, $|\delta_n+2F_n(a^n+b^n)|$ is at least this lower bound minus a term $O(F_n|a|^nn^{-2})$. As $a\neq b$, this term is at most half of the lower bound for large $n$, and the first bound follows.

For the second, we show that $|z^n-1|<2\pi/n$ for infinitely many $n$. By \eqref{eq:shift-cons}, $z=e^{2\pi i\vartheta}$ with $\vartheta=2\tau/\pi$, and by \cite[Theorem~21]{Khi64} there are infinitely many $n\in\N$ with $|n\vartheta-m|<1/n$ for some $m\in\Z$. For these $n$, $|z^n-1|=|e^{2\pi i(n\vartheta-m)}-1|\le2\pi|n\vartheta-m|<2\pi/n$, so \eqref{eq:main3} holds along them, that is,
\[
\delta_n=P_n-\frac{(a-b)^2a^nF_n}{2n^2}+O\bigl(F_n|a|^nn^{-3}\bigr).
\]
As $a\neq b$, the error term is at most $|a-b|^2|a|^nF_n/(4n^2)$ in modulus for large $n$, and the second bound follows.

(b) If $a=b$, then $P_n=-4a^nF_n$, and we may take $q=a$ and $\tau=0$ in \eqref{eq:shift}. Then $v=u$, so $s_v=S$ and hence $\mu_n=\ls$. Thus \eqref{eq:delta-mid} reduces to $\delta_n=-\frac12\gamma_n$, and Proposition~\ref{lem:gap} with $q=a$ gives (b) in the form
\[
\delta_n=-4a^nF_n\Bigl(1-\frac{a^2}{4n^3}+O(n^{-4})\Bigr)=P_n+\frac{a^{n+2}F_n}{n^3}+O\bigl(F_n|a|^nn^{-4}\bigr).
\]

The last statement follows, since (b) gives $\delta_n-P_n=O(F_n|a|^nn^{-3})$ if $a=b$, and the second bound in (a) excludes $\delta_n-P_n=o(F_n|a|^nn^{-2})$ if $a\neq b$.
\end{proof}

\section{Leading terms of the deviations}\label{sec:consequences}

In this section $a,b\in\C\setminus\{0\}$ are fixed with $|a|=|b|$. Proposition~\ref{prop:position} and Theorems~\ref{thm:main} and~\ref{thm:refined} bound the errors in absolute terms. Here we determine when these errors are small compared with the terms they accompany, and so obtain the leading terms of the deviations.

We begin with the first term of \eqref{eq:main2}, that is, with the right-hand side of the Djakov--Mityagin formula \eqref{eq:DM}. Its modulus is $2F_n|a|^n\rho_n$, where $\rho_n=|a^n+b^n|/|a|^n=|1+(a/b)^n|$, while the second term is $O(F_n|a|^n/n)$, so $n\rho_n$ compares the first term with this scale. Let $N_*=\{n\in\N:\ a^n+b^n\neq0\}$. This set is infinite, since it contains $n$ or $2n$ for every $n$. Indeed, $a^n+b^n=0$ gives $(a/b)^{2n}=1$, so $a^{2n}+b^{2n}=2b^{2n}\neq0$.

\needspace{4\baselineskip}
\begin{theorem}[validity of the Djakov--Mityagin formula]\label{thm:validity}
\textup{(a)} If $N'\subseteq N_*$ is infinite, then $\delta_n\sim-2F_n(a^n+b^n)$ as $n\to\infty$ along $N'$ if and only if $n\rho_n\to\infty$ along $N'$.

\noindent\textup{(b)} We have $\delta_n\sim-2F_n(a^n+b^n)$ as $n\to\infty$ along $N_*$ if and only if $a/b$ is a root of unity.

\end{theorem}

\begin{proof}
(a) As $|a|=|b|$, the parallelogram law gives $|a^n-b^n|^2=4|a|^{2n}-|a^n+b^n|^2$, so $|a^n-b^n|/|a|^n=\sqrt{4-\rho_n^2}$. For $n\in N_*$ we have $|a^n+b^n|=|a|^n\rho_n$, so \eqref{eq:main2} and the reverse triangle inequality give
\begin{equation}\label{eq:rel-err}
n\rho_n\,\Bigl|\frac{\delta_n}{-2F_n(a^n+b^n)}-1\Bigr|=\frac{n\,|\delta_n+2F_n(a^n+b^n)|}{2F_n|a|^n}=\frac{|a-b|}{2}\sqrt{4-\rho_n^2}+O(n^{-1}).
\end{equation}
On the left-hand side, $n\rho_n$ multiplies the relative error of $-2F_n(a^n+b^n)$ as an approximation of $\delta_n$, and the right-hand side is at most $|a-b|+O(n^{-1})$. If $n\rho_n\to\infty$ along $N'$, the relative error therefore tends to $0$ along $N'$, that is, $\delta_n\sim-2F_n(a^n+b^n)$ along $N'$.

Conversely, suppose that $n\rho_n$ does not tend to infinity along $N'$. Then there are $B>0$ and infinitely many $n\in N'$ with $n\rho_n\le B$. Along these indices $\rho_n\le B/n\to0$. Hence $a\neq b$, as $\rho_n=2$ if $a=b$, and the right-hand side of \eqref{eq:rel-err} tends to $|a-b|>0$, so it is at least $|a-b|/2$ for large $n$. Dividing \eqref{eq:rel-err} by $n\rho_n>0$ gives $|\delta_n/(-2F_n(a^n+b^n))-1|\ge|a-b|/(2n\rho_n)\ge|a-b|/(2B)$. Thus $\delta_n\sim-2F_n(a^n+b^n)$ fails along $N'$.

(b) Suppose first that $(a/b)^m=1$ for some $m\in\N$. Then $\rho_n$ is $m$-periodic, so it takes only finitely many values. As $\rho_n>0$ exactly for $n\in N_*$, its values on $N_*$ have a minimum $\rho_*>0$, that is,
\begin{equation}\label{eq:rho-lower}
\rho_n\ge\rho_*\qquad\text{for }n\in N_*.
\end{equation}
Hence $n\rho_n\to\infty$ along $N_*$, and part~(a), with $N'=N_*$, gives $\delta_n\sim-2F_n(a^n+b^n)$ along $N_*$.

Suppose now that $a/b$ is not a root of unity. Then no power of $a/b$ equals $-1$, since the square of such a power would be $1$, so $N_*=\N$. Here $a/b=e^{2\pi i\vartheta}$ with $\vartheta=2\tau/\pi$ by \eqref{eq:shift-cons}, and $\vartheta$ is irrational. By Chebyshev's theorem on inhomogeneous approximation \cite[Theorem~24]{Khi64}, infinitely many $n\in\N$ satisfy $|n\vartheta-j-\tfrac12|<3/n$ for some $j\in\Z$. Since $-1=e^{2\pi i(j+\frac12)}$, for these $n$ we have
\[
n\rho_n=n\bigl|e^{2\pi in\vartheta}-e^{2\pi i(j+\frac12)}\bigr|\le2\pi n\bigl|n\vartheta-j-\tfrac12\bigr|<6\pi.
\]
Thus $n\rho_n$ does not tend to infinity along $N_*$, and part~(a) shows that $\delta_n\sim-2F_n(a^n+b^n)$ fails along $N_*$.
\end{proof}

We now compare $\delta_n$ with the sum $P_n$ of the first two terms in \eqref{eq:main2}. The next two propositions assume that $b\neq\bar a$ or $a/b$ is a root of unity. If $a/b$ is not a root of unity, this gives $\operatorname{Re}(a-b)\neq0$. Indeed, then $b\neq a$ and $b\neq\bar a$, while $\operatorname{Re}a=\operatorname{Re}b$ and $|a|=|b|$ would give $|\Img a|=|\Img b|$, that is, $b=a$ or $b=\bar a$.

\begin{proposition}[midpoint deviations]\label{prop:midpoint}
If $b\neq\bar a$ or $a/b$ is a root of unity, then $\delta_n=P_n\bigl(1+O(n^{-1})\bigr)$ as $n\to\infty$, with $P_n\neq0$ for all large $n$. As $n\to\infty$ along the indices with $a^n+b^n=0$,
\begin{equation}\label{eq:cancel}
\delta_n=\frac{2(a-b)\,a^nF_n}{n}\bigl(1+O(n^{-1})\bigr).
\end{equation}
\end{proposition}

\begin{proof}
By \eqref{eq:main2}, $\delta_n-P_n=O(F_n|a|^nn^{-2})$. For the first assertion, it therefore suffices to find $\kappa>0$ such that
\begin{equation}\label{eq:Pn-lower}
|P_n|\ge\frac{\kappa F_n|a|^n}{n}\qquad\text{for large }n.
\end{equation}
Indeed, \eqref{eq:Pn-lower} gives $P_n\neq0$ and $|\delta_n/P_n-1|=|\delta_n-P_n|/|P_n|=O(n^{-1})$ for large $n$.

Suppose first that $a/b$ is a root of unity. For $n\in N_*$, the reverse triangle inequality and \eqref{eq:rho-lower} yield $|P_n|\ge2F_n|a^n+b^n|-(F_n/n)|a-b|\,|a^n-b^n|\ge F_n|a|^n(2\rho_*-2|a-b|/n)$, and the right-hand side is at least $\rho_*F_n|a|^n$ for large $n$. If $n\notin N_*$, then $b^n=-a^n$, so $a\neq b$ and $P_n=2(a-b)a^nF_n/n$. Hence \eqref{eq:Pn-lower} holds with $\kappa=\rho_*$ if $N_*=\N$, and with $\kappa=\min\{\rho_*,2|a-b|\}$ otherwise.

Suppose now that $a/b$ is not a root of unity. Then $\operatorname{Re}(a-b)\neq0$. With $\xi_n=(a-b)/(2n)$ as in \eqref{eq:exp-id}, we have $P_n=-2F_n\bigl[(1-\xi_n)a^n+(1+\xi_n)b^n\bigr]$ and $|\xi_n|\le|a|/n$, so $|1\pm\xi_n|\le2$ for $n\ge|a|$. As $|a^n|=|b^n|$, the reverse triangle inequality and $|1+\xi_n|^2-|1-\xi_n|^2=4\operatorname{Re}\xi_n$ give
\[
\begin{aligned}
|P_n|&\ge2F_n|a|^n\,\bigl||1+\xi_n|-|1-\xi_n|\bigr|\\
&=\frac{8F_n|a|^n\,|\operatorname{Re}\xi_n|}{|1+\xi_n|+|1-\xi_n|}\ge\frac{F_n|a|^n}{n}\,|\operatorname{Re}(a-b)|.
\end{aligned}
\]
Hence \eqref{eq:Pn-lower} holds with $\kappa=|\operatorname{Re}(a-b)|$.

Finally, if $a^n+b^n=0$, then $(a/b)^{2n}=1$, so the first assertion applies, and inserting $P_n=2(a-b)a^nF_n/n$ into it gives the formula \eqref{eq:cancel}.
\end{proof}

In the cases $a=-b$ and $a=ib$ of an example in \cite{DM20}, where $b>0$, the edge formulas show only that $\delta_n=o(|\lambda_n^+-\lambda_n^-|)$ at the indices with $a^n+b^n=0$. Proposition~\ref{prop:midpoint} gives the leading term. For $a=-b$ these indices are the odd $n$, and \eqref{eq:cancel} gives $\delta_n\sim4b^{n+1}F_n/n$. For $a=ib$ they are the $n$ with $n\equiv2\pmod4$, and \eqref{eq:cancel} gives $\delta_n\sim2(1-i)b^{n+1}F_n/n$.

\needspace{4\baselineskip}
We now turn to the two edge deviations $\mu_n-\ls$ and $\mu_n-\lc$, with $q$ as in \eqref{eq:shift}.

\begin{proposition}[edge deviations]\label{prop:edges}
Suppose $a\neq b$.

\noindent\textup{(a)} If $b\neq\bar a$ or $a/b$ is a root of unity, then, as $n\to\infty$,
\[
\mu_n-\ls=\gamma_n\sin^2\Th\bigl(1+O(n^{-1})\bigr),\qquad \mu_n-\lc=-\gamma_n\cos^2\Th\bigl(1+O(n^{-1})\bigr).
\]

\noindent\textup{(b)} As $n\to\infty$ along the indices with $a^n=b^n=q^n$,
\[
\mu_n-\ls=-\frac{(a-b)^2a^nF_n}{2n^2}\bigl(1+O(n^{-1})\bigr),\qquad \mu_n-\lc=-8a^nF_n\bigl(1+O(n^{-2})\bigr),
\]
so $\mu_n$ is nearer to $\ls$ for large $n$. Along the indices with $a^n=b^n=-q^n$, the same holds with $\ls$ and $\lc$ interchanged, and $\mu_n$ is nearer to $\lc$.
\end{proposition}

\begin{proof}
(a) By Proposition~\ref{prop:position}, $\mu_n-\ls=w_n\gamma_n$ and $\mu_n-\lc=-(1-w_n)\gamma_n$, with $w_n$ and $1-w_n$ given by \eqref{eq:weight}. It therefore suffices to find $c>0$ such that
\begin{equation}\label{eq:trig-lower}
|\sin\Th|\ge\frac{c}{n},\qquad|\cos\Th|\ge\frac{c}{n}\qquad\text{for large }n.
\end{equation}
Indeed, by \eqref{eq:trig-lower} an error $O(n^{-2})$ divided by $\sin\Th$ is $O(n^{-1})$, so $\sin\Th+O(n^{-2})=\sin\Th\bigl(1+O(n^{-1})\bigr)$, and the first formula in \eqref{eq:weight} becomes
\[
w_n=\sin^2\Th\bigl(1+O(n^{-1})\bigr)^2\bigl(1+O(n^{-2})\bigr)=\sin^2\Th\bigl(1+O(n^{-1})\bigr),
\]
and the second gives $1-w_n=\cos^2\Th\bigl(1+O(n^{-1})\bigr)$ in the same way.

Suppose first that $\operatorname{Re}(a-b)\neq0$. The identities
\[
|\sin\Th|^2=\sin^2(\operatorname{Re}\Th)+\sinh^2(\Img\Th),\qquad|\cos\Th|^2=\cos^2(\operatorname{Re}\Th)+\sinh^2(\Img\Th)
\]
show that both $|\sin\Th|$ and $|\cos\Th|$ are at least $|\sinh\Img\Th|\ge|\Img\Th|$. As $\Img\Th=\operatorname{Re}(a-b)/(4n)$ by \eqref{eq:theta}, the bound \eqref{eq:trig-lower} holds with $c=|\operatorname{Re}(a-b)|/4$.

In the remaining case, $\operatorname{Re}(a-b)=0$ and $a/b$ is a root of unity. Then $\Th$ is real by \eqref{eq:theta}. If $(a/b)^m=1$ with $m\in\N$, then $(e^{i\tau})^{4m}=1$ by \eqref{eq:shift-cons}, so $\sin n\tau$ and $\cos n\tau$ take only finitely many values, and there is $c_1>0$ such that each of them is either $0$ or of modulus at least $c_1$.

Consider the sine. As $\Th$ and $n\tau$ are real, $|\sin\Th-\sin n\tau|\le|\Th-n\tau|=|a-b|/(4n)$ by \eqref{eq:theta}, so for large $n$ with $\sin n\tau\neq0$ we have $|\sin\Th|\ge c_1-|a-b|/(4n)\ge c_1/2$. If $\sin n\tau=0$, then $\cos n\tau=\pm1$, and the addition formula for the sine gives
\[
|\sin\Th|=|\sin(\Th-n\tau)|\ge\tfrac12|\Th-n\tau|=\frac{|a-b|}{8n}
\]
for large $n$, since $|\sin t|\ge|t|/2$ for real $t$ with $|t|\le1$. The same argument applies to the cosine, except that $|\cos\Th|=|\sin(\Th-n\tau)|$ when $\cos n\tau=0$. Hence \eqref{eq:trig-lower} holds with $c=|a-b|/8$.

(b) At the indices with $a^n=b^n$ we have $(a/b)^n=1$, so \eqref{eq:main3} applies and gives
\[
\delta_n=-4a^nF_n-\frac{(a-b)^2a^nF_n}{2n^2}+O\bigl(F_n|a|^nn^{-3}\bigr).
\]
Further, \eqref{eq:delta-mid} and $\lc=\ls+\gamma_n$ give the two edge deviations $\mu_n-\ls=\delta_n+\frac12\gamma_n$ and $\mu_n-\lc=\delta_n-\frac12\gamma_n$.

Since $a^n=b^n$, we have $q^{2n}=a^nb^n=a^{2n}$, so $q^n=\pm a^n$, and \eqref{eq:gap-lead} gives $\frac12\gamma_n=\pm4a^nF_n\bigl(1+O(n^{-3})\bigr)$. With the expansion of $\delta_n$ this yields
\[
\delta_n\pm\tfrac12\gamma_n=-\frac{(a-b)^2a^nF_n}{2n^2}+O\bigl(F_n|a|^nn^{-3}\bigr),\qquad\delta_n\mp\tfrac12\gamma_n=-8a^nF_n+O\bigl(F_n|a|^nn^{-2}\bigr),
\]
with the upper signs if $q^n=a^n$ and the lower signs if $q^n=-a^n$. If $q^n=a^n$, the left-hand sides are $\mu_n-\ls$ and $\mu_n-\lc$, and factoring out the leading terms, with $a\neq b$, gives the two formulas of the statement. The first deviation is $O(n^{-2})$ times the second, so $\mu_n$ is nearer to $\ls$ for large $n$. If $q^n=-a^n$, the left-hand sides are $\mu_n-\lc$ and $\mu_n-\ls$, and the same holds with $\ls$ and $\lc$ interchanged.
\end{proof}

In the cases $a=-b$ and $a=ib$ of the example in \cite{DM20}, where $b>0$, their formula for $\mu_n-\lambda_n^-$ is not stated at the indices with $a^n=b^n$, where its leading expression vanishes. They give $\mu_n-\lambda_n^+\sim-8b^nF_n$ there. By Proposition~\ref{prop:edges}(b), this is the deviation from the farther edge, so $\lambda_n^-$ is the nearer edge in their labelling for large $n$. For any nonzero $b$, the same proposition shows that the deviation of $\mu_n$ from the nearer edge is asymptotic to $-2b^{n+2}F_n/n^2$ for $a=-b$ and even $n$, and to $i\,b^{n+2}F_n/n^2$ for $a=ib$ and $n\equiv0\pmod4$.

The orders in Propositions~\ref{prop:midpoint} and~\ref{prop:edges}(b) share a source in the phase correction $\Th-n\tau=i(a-b)/(4n)$ of \eqref{eq:theta}. When $a^n+b^n=0$, we have $\cos2n\tau=0$ by \eqref{eq:shift-cons}, so $\cos2\Th$ in \eqref{eq:position} is linear in this correction to leading order, and $\delta_n$ is of order $n^{-1}$ times the gap. When $a^n=b^n$ and $a\neq b$, one of $\sin n\tau$ and $\cos n\tau$ vanishes, so the corresponding weight in \eqref{eq:weight} is the square of the correction to leading order, and the deviation from the nearer edge is of order $n^{-2}$ times the gap.

\section*{Declaration of Generative AI and AI-assisted technologies}

ChatGPT/Codex (OpenAI) was used for language editing. The author takes responsibility for the publication.

\end{document}